\documentclass{amsart}
\newtheorem{lem}{Lemma}[section]

\newtheorem{thm}{Theorem}[section]
\newtheorem{cor}{Corollary}[section]
\newtheorem{defn}{Definition}[section]

\newcommand{\mb}{\mathbb}
\newcommand{\commentout}[1]{}

\newcommand{\mc}{\mathcal}

\newcommand{\Hom}{\rm Hom_{\mathbb R}}

\begin{document}
\title{On Matrix Valued Square Integrable Positive Definite Functions}
\author{Hongyu He}
\address{Department of Mathematics, Louisiana State University, Baton Rouge, LA 70803}
\email{livingstone@alum.mit.edu}
\date{}
\footnote{This research is partially supported by NSF grant DMS-0700809.}
\begin{abstract}
In this paper, we study matrix valued
positive definite functions on a unimodular group. We generalize two important results of Godement on $L^2$ positive definite functions. We show that a matrix-valued {\it continuous} $L^2$ positive definite function can always be written as the convolution of an matrix-valued $L^2$ positive definite function with itself. We also prove that,
given two $L^2$ matrix valued positive definite functions $\Phi$
and $\Psi$, $\int_G Tr(\Phi(g) \overline{\Psi(g)}^t) d g \geq 0$. In addition this integral equals zero if and only if $\Phi * \Psi=0$. Our proofs are operator-theoretic and independent of the group.
\end{abstract}
\maketitle

\section{Introduction}
About 60 years ago, Godement published a paper on square integrable positive definite functions on a locally compact group (~\cite{go}). In his paper, Godement proved that every {\it continuous} square integrable positive definite function has an $L^2$-positive definite square root. He also proved, among others, that the inner product between two positive definite $L^2$-functions must be nonnegative. Godement's results and proofs were quite elegant. 
The purpose of this paper is to extend Godement's theorem to matrix-valued positive definite functions on unimodular groups. Obviously, the diagonal of  matrix-valued positive definite functions must all be positive definite. Yet, there is not much to say about the off-diagonal entries and their relationship with diagonal entries. So Godement's results do not carry easily to the matrix-valued case. In this paper, we generalize Godement's theorems  to matrix-valued positive definite functions (\cite{go} and Ch 13.\cite{dix}). Our results, we believe, are new.\\
\\
Let $M_n(\mathbb C)$ be the set of $n \times n$ matrices. For a matrix $A$, let $[A]_{ij}$ be the $(i, j)$-th entry of $A$. Let $G$
be a unimodular group. A continuous function $\Phi: G \rightarrow M_n(\mathbb C)$ is
said to be {\it positive definite} if for any $\{ \mathbf C_i \in \mathbb
C^n \}_{i=1}^{l}$ and $\{x_i \in G \}_{i=1}^{l}$,
$$\sum_{i,j=1}^{l} (\mathbf C_i)^t \Phi(  x_i^{-1} x_j) \overline{\mathbf
C_j} \geq 0.$$
Take $x_1=e$ and $x_2=g$. The above inequality implies that $\Phi(g)=\overline{\Phi(g^{-1})}^t$ (See for example, Prop. 2.4.6 \cite{kr}). When $n=1$, our definition agrees with the definition of continuous positive
definite functions. We denote the set of continuous matrix-valued positive definite functions by $\mathcal P(G, M_n)$. 
\begin{defn} Let $L_{loc}^1 (G, M_n)$ be the set of $M_n$-valued locally integrable functions on $G$. Let $\Phi \in L_{loc}^1(G, M_n)$ act on $u \in C_c(G, \mathbb C^n)$ by 
$[\lambda(\Phi) (u)(x)]_i =\sum_{j=1}^n \int [\Phi(g)]_{ij} [u(g^{-1} x)]_j d g.$
We write $\lambda(\Phi)(u)(x)=\int \Phi(g) u(g^{-1} x ) d g$.
Clearly, $\lambda(\Phi) u$ is continuous.
We say that $\Phi $ is positive definite if 
$$\langle \lambda(\Phi)(u), u \rangle=\sum_{i=1}^n \int_{G} [\lambda(\Phi) u(g)]_i \overline{[u(g)]_i} d g \geq 0$$
for all $ u \in C_c(G, \mathbb C^n)$.
\end{defn}
We denote the set of matrix-valued positive definite functions by
$\mathbf P(G, M_n)$. Clearly, $\mathbf P(G, M_n) \supset \mathcal P(G, M_n)$ (see Prop 13.4.4 \cite{dix}). 
\begin{defn} A matrix-valued function $\Phi(x)$ is said to be
square integrable, or simply $L^2$ if $[\Phi]_{ij}$ is in $L^2(G)$
for all $(i,j)$. We denote the set of matrix-valued square
integrable function by $L^2(G, M_n)$. Define
$$\langle \Phi, \Psi \rangle= \int_{G} Tr \Phi \overline{\Psi}^t d g .$$
Put $\mathcal P^2(G, M_n)=L^2(G, M_n) \cap \mathcal P(G, M_n)$ and $\mathbf P^2(G, M_n)=L^2(G, M_n) \cap \mathbf P(G, M_n)$.
\end{defn}
Let $\Phi, \Psi \in L_{loc}^1(G, M_n)$. Define the convolution
$$[\Phi * \Psi]_{ij}= \sum_{k=1}^n [\Phi]_{i k} * [\Psi]_{k j}$$
whenever the right hand side is well-defined, i.e., the convolution integral converges absolutely.\\
\\
{\bf Theorem [A]}
{\it 
Let $G$ be a unimodular locally compact group. Let $\Phi \in \mc P^2(G, M_n)$. Then there exists a $\Psi \in \mathbf P^2(G, M_n)$ such that $\Phi=\Psi * \Psi$. }\\
\\
{\bf Theorem [B]}
{\it Let $G$ be a unimodular locally compact group. Let $\Phi, \Psi \in \mathbf P^2(G, M_n)$. Then 
$\langle \Phi, \Psi \rangle \geq 0.$
}\\
\\
{\bf Theorem [C]}
{\it Let $G$ be a unimodular locally compact group. Let $\Phi, \Psi \in \mathbf P^2(G, M_n)$. Then $\langle \Phi, \Psi \rangle=0$ if and only of $\Phi * \Psi =0$. }\\
\\
Our motivation comes from the theory of unitary representations of Lie groups. There are representations that appear as a space of \lq\lq invariant distributions \rq\rq in a unitary representation $(\pi, \mathcal H_{\pi})$. To construct a Hilbert inner product for the invariant distributions, one is often led to investigate whether 
 \begin{equation}~\label{basic}
 \int_G (\pi(g)u, u) d g \geq 0
 \end{equation}
 when $(\pi(g) u, u)$ is $L^1$, $u \in \mc H_{\pi}$ and $G$ is unimodular (\cite{li2} \cite{unit}). For $G=\mathbb R$, the affirmative answer to this question is a direct consequence of Bochner's theorem, namely, the integral of a $L^1$ positive definite function on $\mathbb R$ is nonnegative.  In his thesis \cite{go}, Godement raised this question for $G$ unimodular. It is known that for $G$ amenable,  the Inequality (\!\!~\ref{basic}) is always true (Prop 18.3.6 \cite{dix}). An amenable group is characterized by the fact that the unitary dual is weakly contained in $L^2(G)$. Therefore for $G$ nilpotent,
the Inequality (\!\!~\ref{basic}) holds. \\
\\
Consider the other extreme, namely, $G$ semisimple and noncompact. The inequality (\ref{basic}) is false in its full generality. Yet, applying the result of this paper, we show that Inequality (\!\!~\ref{basic}) holds if $\mc H_{\pi}$ can be written as a tensor product of two $L^2$-representations of $G$ and $u$ is finite in the tensor decomposition. See Theorem \ref{basic1}.
In Cor. \ref{theta}, we give a result about a certain integral related to Howe's correspondence (\cite{howe}). It is more general than the results given in \cite{unit}. 
\section{Convolution Algebras }
Let $G$ be a unimodular locally compact group.
A matrix-valued function on $G$ is said to be in $L^p$ if each entry is in $L^p(G)$. Let $\Phi \in L^1(G, M_n)$. Define $\| \Phi \|_{L^1}= \sum \| [\Phi]_{ij} \|_{L^1}$. Then $L^1(G, M_n)$ becomes a Banach algebra. We have
$$L^1(G, M_n) * L^1(G, M_n) \subseteq L^1(G, M_n); \qquad C_c(G, M_n) * C_c(G, M_n) \subseteq C_c(G, M_n).$$
$$L^2(G, M_n) * L^2(G, M_n) \subseteq \mc BC(G, M_n).$$
Here $\mc BC(G, M_n)$ is the space of bounded continuous functions.\\
\\
For each $u, v \in L^2(G, \mathbb C^n)$, define the standard inner product $\langle u, v \rangle=\int_G \sum_i [u(g)]_i \overline{[v(g)]_i} d g$.
Obviously $L^1(G, M_n)$ acts on $L^2(G, \mathbb C^n)$ by $\lambda$. Then map
$$\lambda: L^1(G, M_n) \rightarrow \mc B(L^2(G, \mathbb C^n))$$
defines a bounded Banach algebra isomorphism. Notice that if $\Phi \in L^2(G, M_n)$ and $u \in L^1(G, \mathbb C^n)$, then $\lambda(\Phi) u \in L^2(G, \mathbb C^n)$. \\
\\
We define the $*$ operation on $L_{loc}(G, M_n)$ by letting
$$[\Phi^*(g)]_{ij}=\overline{[\Phi(g^{-1})]_{j i}}.$$
For $u, v \in C_c(G, \mathbb C^n)$, we have
$\langle \lambda(\Phi) u, v \rangle=\langle u, \lambda(\Phi^*) v \rangle$.
If $\Phi$ is positive definite in $L_{loc}(G, M_n)$, then
$$\langle \lambda(\Phi) u, u \rangle=\overline{\langle \lambda(\Phi) u, u \rangle}=\langle u, \lambda(\Phi)u \rangle=\langle \lambda(\Phi^*) u, u \rangle.$$
Hence $\langle \lambda(\Phi-\Phi^*) u, u \rangle=0$. Let $u= v + t w \, (t \in \mathbb R)$. Then
$$0=\langle \lambda(\Phi-\Phi^*)(v+ t w), v+ t w \rangle= \langle \lambda(\Phi-\Phi^*) v,  t w \rangle+\langle
\lambda(\Phi-\Phi^*) t w, v \rangle.$$
We obtain
 $\langle \lambda(\Phi-\Phi^*) v , w \rangle=\langle \lambda(-\Phi+\Phi^*) w, v \rangle=\langle w, \lambda(-\Phi^*+ \Phi) v \rangle= \overline{\langle \lambda(\Phi-\Phi^*) v ,w \rangle}.$ Hence
 $\langle \lambda(\Phi-\Phi^*) v , w \rangle$ must be real for all $v, w \in C_c(G, \mathbb C^n)$.
 It follows that $\Phi^*=\Phi$. We have\\
 
\begin{lem}~\label{1} Let $\Psi, \Phi \in \mathbf P^2(G, M_n)$. Then $\Phi=\Phi^*$, $\Psi=\Psi^*$ and
$\langle \Phi, \Psi \rangle= Tr(\Phi*\Psi(e))$.
\end{lem}
Proof: The last statement follows from
\begin{equation}
Tr(\Phi*\Psi(e))=\sum_{i, j} \int_{G} [\Phi]_{i j}(g) [\Psi]_{j i}(g^{-1}) d g = \sum_{i, j} \int_{G} [\Phi]_{i j}(g) \overline{[\Psi]_{i j}(g)} d g =\int_{G} Tr \Phi(g) \overline{\Psi(g)^t} d g=\langle \Phi, \Psi \rangle .
\end{equation}

\begin{defn}[Ch 13. \cite{dix}] Let $\Phi \in L^2(G, M_n)$. We say that $\Phi$ is moderated if $\lambda(\Phi)$ on $C_c(G, \mathbb C^n)$ is a bounded operator in the $L^2$ norm, i.e., there is a $M$ such that
$$\|\lambda(\Phi) u \| \leq M \| u \|$$
for any $u \in C_c(G, \mathbb C^n)$. 
\end{defn}

When $\Phi$ is moderated, $\lambda(\Phi)|_{C_c(G, \mathbb C^n)}$ can be extended to a bounded operator on $L^2(G, \mathbb C^n)$  which coincides with the operator $\lambda(\Phi)|_{L^2(G, \mathbb C^n)}$. To see this, let $u_i \rightarrow u$ under the $L^2$-norm with $u_i \in C_c(G, \mathbb C^n)$. Since $\lambda(\Phi)$ is bounded on $C_c(G, M_n)$, $\{ \lambda(\Phi) u_i \}_{i=1}^{\infty}$ yields a Cauchy sequence in $L^2(G, M_n)$. Therefore $\lambda(\Phi) u_i$ converges to an $L^2$-function $v \in L^2(G, \mathbb C^n)$ . In particular, $(\lambda(\Phi) u_i)(g)$ converges in $L^2$-norm to $v(g)$ on any compact subset $K$. On the other hand,  $(\lambda(\Phi) u_i)(g)$ converges to $\lambda(\Phi) u(g)$ uniformly on $G$, in particular, on $K$. Hence $(\lambda(\Phi)u) (g)=v(g)$ for $g \in K$ almost everywhere. It follows that $v(g)=\lambda(\Phi)u (g)$ almost everywhere. Therefore  $\lambda(\Phi) u \in L^2(G, \mathbb C^n)$.  In short, if $\Phi$ is moderated, and $u \in L^2(G, \mathbb C^n)$, then  $\lambda(\Phi)(u) \in L^2(G, \mathbb C^n)$. We retain $\lambda(\Phi)$ to denote the bounded operator on $L^2(G, \mathbb C^n)$. The following is obvious.
\begin{lem}
$\Phi$ is moderated if and only if $[\Phi]_{ij}$ are all moderated in $L^2(G)$.
\end{lem}
Let $M(G, M_n)$ be the space of moderated $L^2$ functions on $G$.
Lemma 13.8.4 ~\cite{dix} asserts that $M(G) * M(G) \subseteq M(G)$ and 
$\lambda|_{M(G)}: M(G) \rightarrow \mc B(L^2(G))$ is an algebra homomorphism.  Therefore,  we obtain
\begin{lem}\label{20} Let $\Phi, \Psi \in L^2(G, M_n)$.
If $\Phi$ and $\Psi$ are moderated, then $\Phi*\Psi$ is also moderated. In addition
$\lambda(\Phi * \Psi)=\lambda(\Phi) \lambda(\Psi).$
\end{lem}

\section{Matrix-Valued Positive Definite Functions}
Let $G$ be a unimodular locally compact group.
Let us recall some basic result from \cite{dix}. Let $\Phi, \Psi \in \mathbf P(G, M_n)$. We define an ordering
$\Phi \preceq \Psi$ if $\Psi-\Phi \in \mathbf P(G, M_n) $. An immediate consequence is that $Tr(\Phi)(e) \leq Tr(\Psi)(e)$. Clearly, $\Phi \preceq \Psi$ if and only if
$$\langle \lambda(\Phi) u, u \rangle \leq \langle \lambda(\Psi) u, u \rangle \qquad (\forall \,\, u \in C_c(G, \mathbb C^n)).$$
\\
For two bounded  operators $X$ and $Y$ in $\mathcal B(\mc H)$, we say that $X \preceq Y$ if $Y-X$ is positive (Ch 2.4 ~\cite{kr}). $Y-X$ is positive implies that $Y-X$ is self-adjoint (Prop. 2.4.6 ~\cite{kr}). If $\Phi$ and $\Psi$ are moderated, then $\Phi \preceq \Psi$ if and only if $\lambda(\Phi) \preceq \lambda(\Psi)$.\\

\begin{thm}[ Prop. 16 \cite{go}] \label{2}
Let $\Phi, \Psi$ be two moderated elements in $\mathbf P^2(G, M_n)$. Suppose that $\Phi* \Psi=\Psi * \Phi$. Then
$\langle \Phi, \Psi \rangle \geq 0$. Let 
$$\Phi_1 \preceq \Phi_2 \preceq \ldots \preceq \Phi_n \preceq \ldots$$
be an increasing sequence of moderated positive definite functions in $L^2(G, M_n)$. Suppose that $\Phi_i$ mutually commute. If $\sup_i \| \Phi_i \|_{L^2} < \infty$, then $\Phi= \lim \Phi_i$ exists in $\mathbf P^2(G, M_n)$.
\end{thm}
The $n=1$ case is proved as  Prop. 16 in \cite{go}. See also 13.8.5, 13.8.4 \cite{dix}.  \\
\\
Proof of Theorem \ref{2}:  $\Phi*\Psi=\Psi* \Phi$ implies that $\lambda(\Phi) \lambda(\Psi)=\lambda(\Psi) \lambda(\Phi)$ as bounded  operators on $L^2(G, \mathbb C^n)$. Since $\lambda(\Phi)$ and $\lambda(\Psi)$ are both positive, they must be self-adjoint. Hence $\lambda(\Phi)\lambda(\Psi)$ must be positive and self-adjoint. In other words, $\lambda(\Phi * \Psi)$ is positive on $L^2(G, \mathbb C^n)$. In particular, it is positive  with respect to $C_c(G, \mathbb C^n)$. Hence $\Phi*\Psi$, as a matrix-valued continuous function, is positive definite. $\Phi*\Psi(e)$ must be a positive semi-definite matrix. By Lemma \ref{1}
$\langle \Phi, \Psi \rangle=Tr (\Phi* \Psi)(e) \geq 0.$\\
\\
Let 
$\Phi_1 \preceq \Phi_2 \preceq \ldots \preceq \Phi_n \preceq \ldots$
be an increasing sequence of moderated positive definite functions in $L^2(G, M_n)$.
For $j \geq i$, notice that $\| \Phi_j \|^2= \| \Phi_i \|^2+  \langle \Phi_i, \Phi_j-\Phi_i \rangle+\langle \Phi_j-\Phi_i, \Phi_i \rangle+ \| \Phi_j-\Phi_i \|^2 \geq \| \Phi_i \|^2.$ 
Since $\sup_i \| \Phi_i \|_{L^2} < \infty$, the sequence $\{ \| \Phi_i \| \}$ is an increasing sequence bounded from above. In particular, it is a Cauchy sequence.  Notice that for $j \geq i$
$$ \| \Phi_j -\Phi_i \|^2 =\| \Phi_j \|^2 -\| \Phi_i \|^2-  \langle \Phi_i, \Phi_j-\Phi_i \rangle-\langle \Phi_j-\Phi_i, \Phi_i \rangle \leq \| \Phi_j \|^2 -\| \Phi_i \|^2.$$
This implies that $\{ \Phi_i  \}$ is a Cauchy sequence in $L^2(G, M_n)$. Let $\Phi$ be the $L^2$-limit of $\{ \Phi_i \}$. For every $u=(u_p) \in C_c(G, \mathbb C^n)$, since $\lambda(u_p)$ is a bounded operator on $L^2(G)$, we obtain
$$([\Phi_i]_{p,q} * u_q, u_p) \rightarrow ([\Phi]_{p,q} * u_q, u_p).$$
It follows that 
$$ 0 \leq \langle \lambda(\Phi_i) u, u \rangle \rightarrow \langle \lambda(\Phi) u, u \rangle .$$
Therefore $\Phi$ is positive definite. $\Box$ \\

\begin{thm}[Thm. 17 ~\cite{go}]\label{10}
 Suppose that $\Phi$ is a moderated element in $\mathbf P^2(G, M_n)$ such that $\Phi \preceq \Theta$ with $\Theta$ a continuous positive definite function. Then there is a unique moderated element $\Psi$ in $\mathbf P^2(G, M_n)$ such that $\Phi=\Psi* \Psi$ in $L^2$. In particular, $\Phi$ equals a continuous positive definite function almost everywhere and $\| \Psi \|^2 \leq Tr(\Theta(e))$.
\end{thm}
In particular, if $\Phi$ is continuous and moderated in $\mathbf P^2(G, M_n)$ , its square root $\Psi$ exists and is unique. \\
\\
The $n=1$ case is established by Godement. Our proof follows from the proof of Theorem 13.8.6 in \cite{dix} for the scalar-valued positive definite $L^2$ functions. The original idea of Godement is to construct an increasing sequence of positive definite moderated elements $\Psi_k$ in $L^2(G)$ that approaches the square root. In Dixmier's book, $\Psi_k= \| \lambda(\Phi) \| p_k(\frac{\Phi}{\| \lambda(\Phi) \|})$. Here  $\{ p_k(t) \}$ is an increasing sequence of nonnegative polynomials on $[0,1]$ such that $p_k(t) \rightarrow \sqrt{t}$ on $[0, 1]$ and $p_k(0)=0$. We shall supply a proof of this fact before we carry out the proof of Theorem \ref{10}.
\begin{lem}\label{taylor1}
There exists a sequence of polynomials 
$$0 \leq p_1(t) \leq p_2(t) \leq \ldots \leq p_k(t) \leq \ldots \leq \sqrt{t} \qquad (t \in [0,1]),$$
such that $p_k(t) \rightarrow \sqrt{t}$ uniformly on $[0,1]$.
\end{lem}
Proof: Consider the function $t^{-\frac{1}{2}}$, $(t \in [0,1])$. Let $q_k(t)$ be the $k$-th Taylor polynomial at $t=1$. Clearly
$$q_{k+1}(t)=q_k(t)+\frac{(\frac{1}{2}) (\frac{3}{2}) \ldots (\frac{2k+1}{2})}{(k+1)!}(1-t)^{k+1}.$$  
Let $p_k(t)=t q_k(t)$. Clearly $p_k(t)$ is an increasing sequence of non-negative continuous functions with limit $\sqrt{t}$ over the interval $[0,1]$. By Taylor's theorem, $p_k(t) \rightarrow \sqrt{t}$ uniformly on $[\epsilon,1]$. On $[0, \epsilon]$, $\sqrt{t}-p_k(t) < \sqrt{t} \leq \sqrt{\epsilon}$. Hence $p_k(t) \rightarrow \sqrt{t}$ uniformly on $[0,1]$.
 $\Box$\\
\\ 
{\bf Proof of Theorem \ref{10}}: Let $\Phi$ be a moderated element in $\mathbf P^2(G, M_n)$ such that $\Phi \preceq \Theta$ with $\Theta$ a continuous positive definite function. Without loss of generality, suppose the operator norm $\| \lambda(\Phi)\|=1$. For any polynomial $p(t)=\sum_{i=0}^{r} a_i t^i$, define
$$p(\Phi)=\sum_{i=0}^r a_i \overbrace{\Phi* \Phi* \ldots * \Phi}^i.$$
Let $\Psi_k=  p_k(\Phi)$ with $p_k(t)$ defined in the last lemma.  Essentially by functional calculus, we will have 
$$\lambda(\Psi_k) \preceq \lambda(\Psi_{k+1}), \qquad \lambda(\Psi_k) \lambda( \Psi_k) \preceq \lambda(\Phi) .$$ 
It follows that $\Psi_k \preceq \Psi_{k+1}$ and $\Psi_k * \Psi_k \preceq \Phi \preceq \Theta$. By taking the value at $e$, we have    $Tr (\Psi_k * \Psi_k(e)) \leq Tr(\Theta(e))$. By Lemma ~\ref{1}, $\|\Psi_k \|$  is bounded by $\sqrt{Tr \Theta(e)}$. Since $\{ \Psi_k \}$ mutually commutes and is an increasing sequence, by Theorem \ref{2}, the $L^2$-limit of $\Psi_k$ exists. Put $\Psi=\lim_{k \rightarrow \infty} \Psi_k$.  By Theorem \ref{2}, $\Psi \in \mathbf P^2(G, M_n)$. We have 
 $$\Psi * \Psi(g)= \lim_{k \rightarrow \infty} \Psi_k * \Psi_k(g),$$
pointwise. Since
$\lim_{k \rightarrow \infty} \lambda(\Psi_k *\Psi_k -\Phi)= 0$ in the operator norm, for $u \in C_c(G, \mathbb C^n)$, we have
$$\lim_{k \rightarrow \infty} \lambda(\Psi_k * \Psi_k ) u = \lambda(\Phi) u$$
in $L^2(G, \mathbb C^n)$. However, the pointwise limit of the left hand side is obviously $\lambda(\Psi * \Psi) u$. Hence $\lambda(\Psi*\Psi-\Phi)u=0$ for every $u \in C_c(G, \mathbb C^n)$. Hence $\Psi* \Psi(g)=\Phi(g)$ almost everywhere. In particular, $\Phi(g)$ is equal to a continuous positive definite function almost everywhere.  \\
\\
Now $\lambda(\Phi)=\lambda(\Psi)^2$ on $C_c(G, \mathbb C^n)$. For any $u \in C_c(G, \mathbb C^n)$, we have
$$\| \lambda(\Psi) u \|^2= \langle \lambda(\Psi)u, \lambda(\Psi) u \rangle= \langle \lambda(\Psi)^2 u, u \rangle= \langle \lambda(\Phi) u, u \rangle \leq \| \lambda(\Phi) \| \| u \|^2.$$
Hence $\lambda(\Psi)$ is bounded on $C_c(G, \mathbb C^n)$ and the function $\Psi$ is moderated. 
By Lemma \ref{20}, $\lambda(\Phi)=\lambda(\Psi)^2$, as bounded  self-adjoint operators on $L^2(G, \mathbb C^n)$. Since $\lambda(\Psi)$ is positive, $\lambda(\Psi)$ is unique as a bounded operator on $L^2(G, \mathbb C^n)$. In particular, $\lambda(\Psi)$ is uniquely defined on $C_c(G, \mathbb C^n)$. Then $\Psi$ must be unique. $\Box$

\section{Square Roots: Proof of Theorem A}
Let $G$ be a unimodular locally compact group. Let $\Phi \in \mc P^2(G, M_n)$.
Now we would like to give a proof of Theorem A. Our proof is somewhat different from the proof of Theorem 13.8.6 given in ~\cite{dix}. The basic idea is the same, namely, to construct a sequence of moderated continuous positive definite functions $\Phi_k \rightarrow \Phi$. Let $\Psi_k$ be the square root of $\Phi_k$. Then the square root of $\Phi$ can be obtained as the $L^2$-limit of $\Psi_k$. The construction is canonical. In our proof, the continuity of $\Phi_k$ is given by Theorem ~\ref{10}. We do not use Cor. 13.7.11 in ~\cite{dix} which requires several more pages of argument. We also wish to point out a major difference. 
In the scalar case $\lambda(\Phi_k)$ acts on $L^2(G)$ and in our case $\lambda(\Phi_k)$ acts on $L^2(G, \mathbb C^n)$ not on
$L^2(G, M_n)$. \\
\\
{\bf Proof of Theorem [A]}: Let $x \in G$.  Let $\rho(x)$ act on $L^2(G, \mathbb C^n)$ by
$(\rho(x) u)(g)= u( g x).$ The action $\rho$ is simply the right regular action. Hence $\rho(x)$ is a unitary operator on $L^2(G, \mathbb C^n)$. If $\Phi \in L_{loc}(G, M_n)$, then obviously
\begin{equation}~\label{inv}
\rho(x) \lambda(\Phi) \rho(x^{-1}) = \lambda(\Phi)
\end{equation}
on $C_c(G, \mathbb C^n)$. \\
\\
 Let $\Phi \in \mathcal P^2(G, M_n)$. Then $\lambda(\Phi)|_{C_c(G, \mathbb C^n)}$ is a positive symmetric operator densely define on $L^2(G, \mathbb C^n)$, by the definition of positive definiteness of $\Phi$. Let $\Lambda(\Phi)$ be the Friedrichs extension of $\lambda(\Phi)|_{C_c(G, \mathbb C^n)}$. Then $\Lambda(\Phi)$ is an (unbounded) positive and self-adjoint operator (Ch. 5.6. ~\cite{kr}).
 By Equation (\!\!~\ref{inv}), we must have $\rho(x) \Lambda(\Phi)\rho(x^{-1})= \Lambda(\Phi)$. \\
\\
Let
$\Lambda(\Phi)=\int_0^{\infty} t d P$ be the spectral decomposition. Here $P$ is a projection-valued measure on the Borel subsets of $\mathbb R$. In other words, for every $B$ a Borel subset of $\mathbb R$, there is a projection $P(B)$ on $L^2(G, \mathbb C^n)$. Then $\rho(x) \Lambda( \Phi) \rho(x^{-1})=\int_0^{\infty} t d  [\rho(x) P \rho(x^{-1})]$. Notice here that $\rho(x)$ is unitary. Hence $\rho(x) P(B) \rho(x^{-1})$ remains a projection. The uniqueness of the spectral decomposition of self-adjoint operators implies that $\rho(x) P(B) \rho(x^{-1})=P(B)$. Since $P(B)$ is bounded, we have $\rho(x) P(B)= P(B) \rho(x)$ for any Borel subset $B$ and for any $x \in G$. \\
\\
Let $[\Phi]_{*j}$ be the $j$-th column vector of $\Phi$. Fix a Borel subset $B$.
Define
$\Phi_B$ by letting the $j$-th column vector to be 
$[\Phi_B]_{* j} = P(B) [\Phi]_{* j}.$ Clearly $\Phi_B \in L^2(G, M_n)$.\\
\\
{\bf Claim 1: $\lambda(\Phi_B)= P(B) \lambda(\Phi) $ on $C_c(G, \mathbb C^n)$}. \\
\\
Proof: Let $u \in C_c(G, \mathbb C^n)$. Then 
$$[\lambda(\Phi) u ](g)=\sum_j \int_{x \in G} [\Phi]_{* j}(g x^{-1}) [u]_j(x) d x=\sum_j \int_{x \in G}(\rho(x^{-1}) [\Phi]_{* j})(g) [u]_j(x) d x.$$
For any $x \in G$,  we have
$$(P(B) \rho(x^{-1}) [\Phi]_{* j})(g)=(\rho(x^{-1}) P(B) [\Phi]_{* j})(g)=
(P(B) [\Phi]_{* j})(g x^{-1}).$$
Since $P(B)$ is a bounded operator on $L^2(G, \mathbb C^n)$, $[\Phi]_{*j} \in L^2(G, \mathbb C^n)$ and  $[u]_j(x) \in L^1(G)$, we have
\begin{equation}
\begin{split}
[P(B) (\lambda(\Phi) u) ](g)= &  P(B)  \int \sum_j  (\rho(x^{-1}) [\Phi]_{* j})(g) [u]_j(x) d x \\
= &  \int \sum_{j} (P(B) \rho(x^{-1})[\Phi]_{* j})(g) [u]_j(x) d x \\
= &  \int \sum_{j} (\rho(x^{-1})P(B) [\Phi]_{* j})(g) [u]_j(x) d x \\
= &  \int \sum_{j} (P(B) [\Phi]_{* j})(g x^{-1}) [u]_j(x) d x \\
= & \int \sum_j ([\Phi_B]_{* j})(g x^{-1}) [u]_j(x) d x \\
= & \int \Phi_B(g x^{-1}) u(x) d x \\
= & (\lambda(\Phi_B) u) (g)
\end{split}
\end{equation}
Our claim is proved. \\
\\
Observe that $P(B) \lambda(\Phi)= P(B) \Lambda(\Phi)$ on $C_c(G, \mathbb C^n)$ and $P(B) \Lambda(\Phi)$ is positive and bounded. Therefore $\lambda(\Phi_B)=P(B) \lambda(\Phi)$ is bounded on  $C_c(G, \mathbb C^n)$ and positive with respect to $C_c(G, \mathbb C^n)$. Hence $\Phi_B$ is moderated and positive definite. We must have
$\lambda(\Phi_B)= P(B) \Lambda(\Phi)$ on $L^2(G, \mathbb C^n)$. In addition if $B_1 \supset B_2$ 
$$\lambda(\Phi_{B_1}-\Phi_{B_2})=(P(B_1)-P(B_2)) \Lambda(\Phi)$$
on $C_c(G, \mathbb C^n)$ and the right hand side is positive and self adjoint. Hence
$\Phi_{B_1} \succeq \Phi_{B_2}$. Similarly $\Phi_{B_1} \preceq \Phi$.\\
\\
For each positive integer $k$, define
$\Phi_k=\Phi_{[0, k]}$. We then  obtain an increasing sequence of moderated positive definite functions
$$ \Phi_1 \preceq \Phi_2 \preceq \ldots \preceq \Phi_k \preceq \ldots (\preceq \Phi) .$$
Due to the way $[\Phi_k]_{* j}$ are defined, $\Phi_k \rightarrow \Phi$ in $L^2$-norm. We have
\begin{lem}\label{p2}
 Let $G$
be a unimodular group.
Every $\Phi \in \mathbf P^2(G, M_n)$ is a $L^2$-limit of an increasing sequence of mutually commutative moderated elements in $\mathbf P^2(G, M_n)$.
\end{lem}
The $n=1$ case was proved by Godement as Prop. 14 in \cite{go}. \\
\\
Since $\Phi_k$ is moderated in $\mathbf P^2(G, M_n)$ with $\Phi_k \preceq \Phi$, by Theorem ~\ref{10}, there is a moderated element $\Psi_k \in \mathbf P^2(G, M_n)$ such that $\Phi_k=\Psi_k * \Psi_k$ almost everywhere.  Without loss of generality, suppose that $\Phi_k=\Psi_k * \Psi_k$ pointwise. Notice that both $\lambda(\Phi_k)$ and $\lambda(\Psi_k)$ can be regarded as positive bounded self-adjoint operators on the Hilbert space  $L^2(G, \mathbb C^n)$. By Lemma \ref{20}, as  bounded self-adjoint operators  on $L^2(G, \mathbb C^n)$,
$\lambda(\Phi_k)=\lambda(\Psi_k)^2$. We have 
\begin{equation}\label{1000}
\lambda(\Psi_k)=\int_0^k \sqrt{t} d P(t).
\end{equation}
In particular, $\lambda(\Psi_k)$ is uniquely defined on $C_c(G, \mathbb C^n)$. Therefore $\Psi_k$ is unique and satisfies Equation \ref{1000}. By functional calculus,   $\{\lambda(\Psi_k) \}$ mutually commute and yield an increasing sequence of positive bounded self-adjoint  operators on $L^2(G, \mathbb C^n)$. Restricted to $C_c(G, \mathbb C^n)$, it is easy to see that $\{\Psi_k \}$ must mutually commute and 
$$\Psi_1 \preceq \Psi_2 \preceq \ldots \preceq \Psi_k \preceq \ldots.$$
\\
Observe that  $\| \Psi_k \|^2= Tr(\Psi_k * \Psi_k(e)) \leq Tr(\Phi(e))$. By Theorem \ref{2}, $\{ \Psi_k \}$ converges in $L^2(G, M_n)$. Let $\Psi_k  \rightarrow \Psi$ in $L^2(G, M_n)$. By Theorem \ref{2}, $\Psi \in \mathbf P^2(G, M_n)$. Notice that $\Psi_k \in L^2(G, M_n)$. Then $\Phi_k=\Psi_k * \Psi_k$ converges uniformly to $\Psi*\Psi$. Since $\Phi_k|_K \rightarrow \Phi|_K$ in $L^2(K, M_n)$ for any compact set $K$, $\Phi|_K=\Psi*\Psi|_K$ almost everywhere. Therefore $\Phi=\Psi*\Psi$ almost everywhere. Since $\Phi$ is continuous, $\Phi=\Psi* \Psi$. Theorem A is proved. $\Box$

\section{Nonnegative Integral: Proof of Theorem B}
 Let $G$
be a unimodular group. Let $\Phi, \Gamma \in \mathbf P^2(G, M_n)$. We want to prove that
$$\langle \Phi, \Gamma \rangle \geq 0.$$
The main idea of the proof here is essentially due to Godement ( Prop.18 ~\cite{go}).  We start with the following lemma. \\

\begin{lem}~\label{p3} Let $G$ be a unimodular group. Every $\Phi$ in $\mathbf P^2(G, M_n)$ is a limit of an increasing sequence of moderated elements in $\mc P^2(G, M_n)$ under the $L^2$ norm.
\end{lem}
Proof: By Lemma \ref{p2}, it suffices to show that every moderated element $\Phi$ in $\mathbf P^2( G, M_n)$  is the $L^2$ limit of an increasing sequence of moderated elements in $\mc P^2(G, M_n)$. Without loss of generality, suppose that $\| \lambda(\Phi) \|=1$. Let $r_k(t)$ be the $k-$th Taylor polynomial of $\frac{1}{t}$ at $t=1$. We define $q_k(t)=t^2 r_k(t)$. Then $q_k(t)$ is an increasing sequence of nonnegative polynomial functions on $[0,1]$ such that $q_k(t) \rightarrow t$ uniformly on $[0,1]$ (c.f. Lemma \ref{taylor1}).\\
\\
\commentout{Let $q_k(t)$ be an increasing sequence of nonnegative polynomial functions on $[0, 1]$ such that $ t^2 | q_k(t)$ and $q_k(t) \rightarrow t $ uniformly on $[0,1]$. $q_k(t)$ can be chosen to be $t^2 r_k(t)$ with $r_k(t)$ the $k$-th Taylor polynomial of $\frac{1}{t}$ at $t=1$.} 
Let $\Phi_k=q_k(\Phi)$. Then $\lambda(\Phi_k)=q_k(\lambda(\Phi))$ is an increasing sequence of positive self-adjoint operators that approaches $\lambda(\Phi)$. Obviously, $\Phi_k(g)$ is positive definite. Since $\lambda(\Phi)$ extends to a bounded operator on $L^2(G, \mb C^n)$, $\lambda(\Phi_k)$ also extends to a bounded operator on $L^2(G, \mb C^n)$. Hence $\Phi_k$ is  moderated. Since $\Phi * \Phi$ is  continuous, $\Phi_k=q_k(\Phi)$ is always continuous.  Therefore, $\{ \Phi_k \}$ is  an increasing sequence of continuous moderated positive definite functions.\\
\\
Notice that $\lambda(\Phi_k)$, $\lambda(\Phi)$ all mutually commute.  Since $\lambda(\Phi_k) \preceq \lambda(\Phi)$, $(\lambda(\Phi_k))^k \preceq (\lambda(\Phi))^k$.  Hence $\Phi_k * \Phi_k \preceq \Phi * \Phi$. This implies $Tr(\Phi_k * \Phi_k (e)) \leq Tr(\Phi*\Phi(e))$. By a similar argument in the proof of Theorem \ref{2},  $\| \Phi_k \| \leq \| \Phi \|$.
By Theorem \ref{2}, let $\Psi$ be the $L^2$-limit of $\Phi_k$. For any  $u \in C_c(G, M_n)$, 
$\lambda(\Psi) u=\lim_{k \rightarrow \infty} \lambda(\Phi_k) u$ pointwise, and $\lim_{k \rightarrow \infty} \lambda(\Phi_k) u =\lambda(\Phi) u$ in $L^2$-norm. It follows that $\Psi=\Phi$ almost everywhere. Therefore $\Phi_k \rightarrow \Phi$ in $L^2$-norm.\\
\\
 We have obtained an increasing sequence of moderated elements in $\mc P^2(G, M_n)$ such that $\Phi_k \rightarrow \Phi$ in $L^2$-norm. $\Box$ \\

\begin{lem} Let $\Phi_1$ be a moderated element in $\mathbf P^2(G, M_n)$ and $\Phi_2 \in \mathcal P^2(G, M_n)$. We have
$$ \langle \Phi_1, \Phi_2 \rangle \geq 0.$$
\end{lem}
Proof: Suppose that $\Phi_2=\Psi * \Psi$ with $\Psi \in \mathbf P^2(G, M_n)$. Then
$$ \langle \Phi_1, \Phi_2 \rangle= Tr(\Phi_1 * \Phi_2(e))=Tr(\Phi_1 * \Psi * \Psi(e))= \langle \lambda(\Phi_1) \Psi,  \Psi \rangle=\sum_{i=1}^n \langle \lambda(\Phi_1) [\Psi]_{ * i}, [\Psi ]_{* i} \rangle.$$
Notice that $\lambda(\Phi_1)$ is a bounded positive self adjoint operator. We have
$ \langle \Phi_1, \Phi_2 \rangle \geq 0.$ $\Box$\\
\\
{\bf  Proof of Theorem B}: For $\Phi, \Gamma \in \mathbf P^2(G, M_n)$, let $\Phi_{\alpha}$ be a sequence of moderated element in $\mathbf P^2(G, M_n)$ with $L^2$-limit $\Phi$ and $\Gamma_{\beta}$ be a sequence of elements in $\mathcal P^2(G, M_n)$ with $L^2$-limit $\Gamma$. Then we have
$$\langle \Phi, \Gamma \rangle= \lim_{\alpha, \beta \rightarrow \infty} \langle \Phi_{\alpha}, \Gamma_{\beta} \rangle  \geq 0.$$
Theorem B is proved.  $\Box$
\section{Zero Integral}
Let $\Phi,\Psi \in \mathbf P^2(G, M_n)$. If $\Phi * \Psi=0$, we have $\langle \Phi, \Psi \rangle= Tr(\Phi*\Psi(e))=0$. Now we would like to show that the converse is also true.
\begin{thm} Let $G$ be a unimodular locally compact  group.
Let $\Phi, \Psi \in \mathbf P^2(G, M_n)$. If $\langle \Phi, \Psi \rangle=0$, then $\Phi * \Psi =0$.
\end{thm}
Proof: By Lemma ~\ref{p2}, let $\Phi_{m}$ be an increasing sequence of moderated elements in $\mathbf P^2(G, M_n)$ such that $\| \Phi_{m}-\Phi \| \leq \frac{1}{m}$. By Lemma ~\ref{p3}, let $\Psi_{p}$ be an increasing sequence in $\mathcal P^2(G, M_n)$ such that $\| \Psi_{p} - \Psi\| \leq \frac{1}{p}$. Then 
$$0=\langle \Phi, \Psi \rangle= \langle \Phi-\Phi_m, \Psi \rangle+ \langle \Phi_m, \Psi \rangle \geq \langle \Phi_m, \Psi \rangle \geq \langle \Phi_m, \Psi_p \rangle \geq 0.$$
Hence all the inequalities here must be equalities.
 Suppose that $\Psi_{p}=\Theta_p * \Theta_p$ with $\Theta_p \in \mathbf P^2(G, M_n)$. Then
$$ 0=\langle \Phi_m, \Psi_p \rangle= Tr(\Phi_m * \Theta_p * \Theta_p(e))=\sum_{i} \langle \lambda(\Phi_m) [\Theta_p]_{* i}, [\Theta_p]_{* i} \rangle.$$
Since $\lambda(\Phi_m)$ is a positive operator on $L^2(G, \mathbb C^n)$, $\langle \lambda(\Phi_m) [\Theta_p]_{* i}, [\Theta_p]_{* i} \rangle =0$. Thus $\lambda(\Phi_m) [\Theta_p]_{* i}=0$ in $L^2(G, \mathbb C^n)$. It follows that $\Phi_m * \Theta_p=0$
in $L^2(G, M_n)$. Since $\Phi_m * \Theta_p(g)$ is a continuous function, $\Phi_m * \Theta_p(g)=0$ for all $g \in G$. Hence $\Phi_m * \Psi_p(g)=\Phi_m * \Theta_p * \Theta_p (g)= 0$. Since $\Phi_m \rightarrow \Phi$ and $\Psi_p \rightarrow \Psi$ in $L^2(G, M_n)$, we have $\Phi_m * \Psi_p(g) \rightarrow \Phi * \Psi(g)$.
 Therefore $ \Phi* \Psi(g)=0$ for all $g$. $\Box$
 
\begin{cor} Let $G$ be a locally compact unimodular group.
Let $\Phi, \Psi \in \mathbf P^2(G, M_n)$. Then $\langle  \Phi, {\Psi} \rangle \geq 0$ and $\langle \Phi, \Psi \rangle=0$ if and only if $\Phi*\Psi=0$.
\end{cor}

\section{Applications in Representation Theory}
Let $G$ be a unimodular group. We call a unitary representation $(\pi, \mc H)$ of $G$ $L^p$ if there is a cyclic vector $u$ in $\mc H$ such that
$(\pi(g)u, u)$ is $L^p$. A $L^p$ unitary representation has a $G$-invariant dense subspace with $L^p$-matrix coefficients.\\

\begin{thm}\label{basic1} Let $G$ be a unimodular locally compact group and $(\pi, \mc H)$ be a unitary representation of $G$. Suppose that $(\pi_1, \mc H_1)$ and $(\pi_2, \mc H_2)$ are two $L^2$-unitary representations of $G$ such that
$$(\pi, \mc H) \cong (\pi_1 \otimes \pi_2, \mc H_1 \hat{\otimes} \mc H_2).$$
 Let $u=\sum_{i=1}^n u_1^{(i)} \otimes u_2^{(i)}$ such that  matrix coefficients with respect to $\{u_1^{(i)}\}$ and $\{ u_2^{(i)} \}$ are all $L^2$. Then
$$\int_{G}(\pi(g) u, u) d g= \sum_{i, j=1}^n \int_G(\pi_1(g) u_1^{(i)}, u_1^{(j)})(\pi_2(g) u_2^{(i)}, u_2^{(j)}) d g \geq 0.$$
\end{thm}
Proof: Observe that $\Phi_1$ defined by $[\Phi_1]_{i j}=(\pi_1(g) u_1^{(i)}, u_1^{(j)})$ is  square integrable and positive definite. Similarly, $\Phi_2 \in L^2(G, M_n)$ defined by $[\Phi_2]_{i j}=( u_2^{(i)}, \pi_2(g) u_2^{(j)})$ is square integrable and positive definite. This theorem follows easily from Theorem B. $\Box$\\
\\
Now we shall apply our result to Howe's correspondence (\cite{howe}). Let $(G(m), G^{\prime}(n))$ be a dual reductive pair in $Sp$. Let $(G^{\prime}(n_1), G^{\prime}(n_2))$ be two $G^{\prime}$-groups diagonally embedded in $G^{\prime}(n)$ with $n_1+n_2=n$. Then $(G(m), G^{\prime}(n_i))$ is a dual reductive pair in some $Sp^{(i)}$ such that $(Sp^{(1)}, Sp^{(2)})$ are diagonally embedded in $Sp$. Let $\omega_i$ be the oscillator representation of
$\widetilde{Sp^{(i)}}$. Let $\omega$ be the oscillator representation of $Sp$. Then $\omega$ can be identified with $\omega_1 \otimes \omega_2$. This identification preserves that actions of $G(m)$ and $G^{\prime}(n_i)$.\\
\\
Now suppose that the matrix coefficients of $\omega_1|_{\tilde{G}(m)}$ with respect to the Schwartz space are $L^2$. Let $\pi$ be an irreducible unitary representation of $\tilde{G}(m)$. Suppose that the matrix coefficients for $\omega_2^{\infty}|_{\tilde {G}(m)} \otimes \pi^{\infty}$ are all square integrable. Then for any $v \in \pi^{\infty}$, $u_1^{(j)} \in \omega_1^{\infty}, u_2^{(j)} \in \omega_2^{\infty}$ with $j \in [1, N]$, we have
\begin{equation}~\label{nonvanish}
\begin{split}
 & \int_{\tilde G(m)} (\omega(g) (\sum u_1^{(j)} \otimes u_2^{(j)}), (\sum u_1^{(k)} \otimes u_2^{(k)}))( \pi(g) v, v) d g \\
=  & \sum_{j, k} \int_{\tilde G(m)} (\omega_1(g) u_1^{(j)},  u_1^{(k)})(\omega_2(g) u_2^{(j)}, u_2^{(k)})(\pi(g) v, v) d g. 
\end{split}
\end{equation}
By Theorem \ref{basic1}, this integral must be nonnegative. \\

\begin{cor}~\label{theta}
Consider a dual reductive pair $(G(m), G^{\prime}(n))$ in $Sp$. Let $n=n_1+n_2$. Let $(G(m), G^{\prime}(n_i))$ be a dual reductive pair in $Sp^{(i)}$. Let $\omega_i$ be the oscillator representation of $Sp^{(i)}$. Let $\pi$ be an irreducible unitary representation of $\tilde{G}(m)$. Suppose that the matrix coefficients with respect to $\omega_1^{\infty}|_{\tilde G(m)}$ and $\omega_2^{\infty}|_{\tilde G(m)} \otimes \pi^{\infty}$ are square integrable.   Let $\xi \in \omega_1^{\infty} \otimes \omega_2^{\infty}$ and $u \in \pi^{\infty}$, then
$$\int_{\tilde{G}(m)} (\omega(g) \xi, \xi)(\pi(g) u, u) d g \geq 0.$$

\end{cor}
This Corollary holds for both p-adic groups and real groups. See \cite{li2}  \cite{theta} \cite{unit} for the importance of this integral in Howe's correspondence (\cite{howe}). In particular, under the hypothesis of the Corollary, Howe's correspondence preserves unitarity.

\end{document}